\documentclass[12pt]{article}
\input epsf.tex

\usepackage{amsmath}
\usepackage{amsthm}
\usepackage{amsfonts}
\usepackage{amssymb}
\usepackage{graphicx}
\usepackage{latexsym}

\usepackage{amsmath,amsthm,amsfonts,amssymb}

\usepackage{epsfig}

\usepackage{amssymb}
\usepackage[all]{xy}
\xyoption{poly}
\usepackage{fancyhdr}
\usepackage{wrapfig}

\theoremstyle{plain}
\newtheorem{thm}{Theorem}[section]
\newtheorem{prop}[thm]{Proposition}
\newtheorem{lem}[thm]{Lemma}
\newtheorem{cor}[thm]{Corollary}

\theoremstyle{definition}
\newtheorem{defn}{Definition}
\theoremstyle{remark}
\newtheorem{remark}{Remark}

\advance \topmargin by -\headheight
\advance \topmargin by -\headsep
\evensidemargin \oddsidemargin
\def\cB{{\mathcal B}}

\def\tB{{\tilde B}}

\def\C{{\mathbb C}}

\def\cc{{\curvearrowright}}

\def\det{{\textnormal{det}}}

\def\hf{{\hat f}}

\def\tG{{\tilde G}}

\def\gr{{\textnormal{Gr}}}

\def\bh{{\bar h}}

\def\N{{\mathbb N}}

\def\tnu{{\tilde \nu}}

\def\per{{\textrm{Fix}}}

\def\tphi{{\tilde{\phi}}}

\def\R{{\mathbb R}}

\def\supp{{\textrm{supp}}}
\def\T{{\mathbb T}}

\def\tr{{\textrm{tr}}}

\def\chix{{\raise.5ex\hbox{$\chi$}}}

\def\cV{{\mathcal V}}

\def\Z{{\mathbb Z}}
\def\bz{{\bf{0}}}

\begin{document}
\title{Entropy for expansive algebraic actions of residually finite groups }
\author{Lewis Bowen\footnote{supported in part by NSF grant DMS-1000104} \\ Texas A\&M University}
\begin{abstract}
We prove a formula for the sofic entropy of expansive principal algebraic actions of residually finite groups, extending recent work of Deninger and Schmidt.
\end{abstract}
\maketitle

\noindent
\section{ Introduction }


Let $\Gamma$ be a countable discrete amenable group. Let $f$ be an element in the group 
ring $\Z\Gamma$. Then $\Gamma$ acts on the abelian group $\Z\Gamma/\Z\Gamma f$ as automorphisms via left 
translation and hence acts on its Pontryagin dual (a compact metrizable abelian 
group) $X_f :=\widehat{\Z\Gamma/\Z\Gamma f}$ by automorphisms. 


The purpose of this paper is to compute the sofic entropy of $X_f$ (with respect to Haar measure) when $\Gamma$ is a residually finite group and the action is expansive (i.e., $f$ is invertible in $l^1(\Gamma)$). Sofic entropy is a generalization of classical entropy to sofic groups introduced in [Bo09b]. The primary novelty is that many non-amenable groups are sofic.

Here is a brief history of related results. The entropy of a single automorphism of a compact separable group was worked out in a series of papers in the 1960's (cf. [LW88,\S 1]) culminating in Yuzvinskii's general result [Yu67]. For $\Z^d$-actions, Lind, Schmidt, and Ward [LSW90] calculated the topological entropy $h(\alpha_f)$ of $X_f$ in terms of the Mahler measure 
of $f$ (for any $f\in \Z\Gamma$), and this is the main step of their calculation for the topological entropy 
of any action of $\Z^n$ on a compact metrizable group as automorphisms. 
In [FK52] Fuglede and Kadison introduced a determinant for invertible elements in 
a unital $C^*$ -algebra $A$ with respect to a tracial state $\tr$. For any invertible $f$ in $A$, 
its Fuglede-Kadison determinant is defined as 
$\det f := \exp(\tr \log 
|f|)$, 
where $|f| = (f^*f)^{1/2}$ 
is the absolute part of $f$. 
 [De06] proved $h(\alpha_f)=\log \det (f)$ when $\Gamma$ is amenable with a strong F\o lner sequence, $f$ is positive in the group von Neumann algebra of $\Gamma$ and is invertible in $l^1(\Gamma)$. [DS07] showed $h(\alpha_f)=\log \det (f)$ when $\Gamma$ is residually finite and amenable and $f$ is invertible in $l^1(\Gamma)$.


In this paper we are concerned exclusively with residually finite groups, which are a special class of sofic groups. Therefore, we give only a special case of the definition of sofic entropy next. We will say that a sequence $\{\Gamma_i\}_{i=1}^\infty$ of finite-index normal subgroups of $\Gamma$ converges to the identity $e$ (written $\lim_{i\to\infty}\Gamma_i = \{e\}$) if  
$$\bigcap_{n \ge 1} \bigcup_{i=n}^\infty \Gamma_i = \{e\}.$$  
$\Gamma $ is residually finite  iff  it has such a sequence. From now on, assume $\Gamma $ is residually finite and $\{\Gamma_i\}_{i=1}^\infty$  is a sequence of finite-index normal subgroups of $\Gamma$ with $\lim_{i\to\infty} \Gamma_i =  \{e\}$. Let $(X,\mu)$  be a standard probability space on which $\Gamma$  acts by measure-preserving transformations.  Let $\phi:X \to A$  be a measurable map into a finite set. In order to define the entropy rate of $\phi$  (with respect to $\{\Gamma_i\}_{i=1}^\infty$ )  we need to compare $\phi$  with measurable maps $\psi:Y \to A$  where $(Y,\nu)$  is any other standard probability space on which $\Gamma$  acts by measure-preserving transformations.  So let $W \subset \Gamma$  be a finite set ($W$ is for {\it window}). Let $A^W$  be the set of functions from $W$ to $A$. Define $\phi^W: X \to A^W$  by $\phi^W(x)(w)=\phi(wx)$. Let $\phi^W_* \mu$  be the pushforward measure on $A^W$.  Define $\psi^W_*\nu$ similarly. The $l^1$-distance between $\phi^W_* \mu$ and $\psi^W_*\nu$ is defined by
$$||\phi^W_* \mu-\psi^W_*\nu||_1:=\sum_{j \in A^W} \big| \phi^W_* \mu(\{j\}) - \psi^W_* \nu(\{j\})\big|.$$

It will be convenient to consider the quotient space $\Gamma_i\backslash \Gamma$ with the left-action of $\Gamma$ given by $\gamma(\Gamma_i\beta):=\Gamma_i \beta \gamma^{-1}$ for any $\gamma \in \Gamma$ and $\Gamma_i\beta \in \Gamma_i\backslash \Gamma$. Let $\zeta_i$ be the uniform probability measure on $\Gamma_i\backslash \Gamma$.

Let $A^{\Gamma_i\backslash\Gamma}$ be the set of functions from $\Gamma_i\backslash\Gamma$ to $A$. Now  define the entropy rate of $ \phi$  with respect to the sequence $\{\Gamma_i\}_{i=1}^\infty$  by
$$h\big(\{\Gamma_i\}_{i=1}^\infty,\phi\big):= \inf_{W \subset \Gamma} \inf_{\epsilon >0} \limsup_{i\to\infty} \frac{\log\Big(\big| \{ \psi\in A^{\Gamma_i \backslash \Gamma}:~ ||\phi^W_* \mu-\psi^W_*\zeta_i||_1<\epsilon\} \big|\Big)}{[\Gamma:\Gamma_i]}$$
 where the first infimum  is over all  finite  $W \subset \Gamma$. 
 
  Recall that  $\phi $ is {\em  generating } if the smallest $\Gamma$-invariant sigma algebra  on $X$  for which $\phi$  is measurable is the sigma algebra of all measurable sets up to sets of measure zero. Part of the main theorem of [Bo09b]  is:
  \begin{thm}\label{thm:K}
  If $\phi:X \to A$  and $ \psi:X \to B$  are  generating where $A$ and $B$  are finite then $h\big(\{\Gamma_i\}_{i=1}^\infty,\phi\big)=h\big(\{\Gamma_i\}_{i=1}^\infty,\psi\big)$.
   \end{thm}
The  number     $h\big(\{\Gamma_i\}_{i=1}^\infty,\phi\big)$  is called the {\em  entropy}  of  the  action $G\cc (X,\mu)$  with respect to $\{\Gamma_i\}_{i=1}^\infty$. It is denoted by $h\big(\{\Gamma_i\}_{i=1}^\infty,X,\mu\big)$. If $\Gamma$  is amenable then this number is the same as the  classical entropy of the action. If $\Gamma$ is non-amenable, this number may depend on the choice of $\{\Gamma_i\}_{i=1}^\infty$.

 The main theorem of this paper is:
 \begin{thm}\label{thm:main}
 If $\Gamma$ is non-amenable and $f \in \Z\Gamma$  is invertible in $l^1(\Gamma)$ then  for any sequence $\{\Gamma_i\}_{i=1}^\infty$ of finite-index normal subgroups of $\Gamma$ with $\lim_{i\to\infty}\Gamma_i=\{e\}$,
 $$h\big(\{\Gamma_i\}_{i=1}^\infty, X_f, \mu_f\big) = \log \det (f)$$
 where $\mu_f$  denotes Haar  measure on $X_f$.
 \end{thm}

In a separate paper, it will be shown that if $\Gamma$ is amenable then for any probability-measure preserving action $\Gamma \cc(X,\mu)$, the classical measure-theoretic entropy rate of a finite observable $\phi:X \to A$, $h(\phi)$ equals $h\big(\{\Gamma_i\}_{i=1}^\infty,\phi\big)$ which also equals $\bh\big(\{\Gamma_i\}_{i=1}^\infty,\phi\big)$ (which is defined in \S \ref{sec:alternative}). Assuming this result, the proof of theorem \ref{thm:main} shows that the conclusion of theorem \ref{thm:main} still holds if $\Gamma$ is amenable. In this way, it is possible to obtain an alternative proof of the main theorem of [DS07]. Note that that paper makes crucial use of tools from topological entropy theory (spanning sets, separating sets, open covers, etc.) that are not available when working with non-amenable groups. 

    The proof of theorem \ref{thm:main} shows more. For any subgroup $\Gamma' <\Gamma$,  let $\per(\Gamma',X_f)$ be the set of all $x\in X_f$ such that $\gamma x = x$ for all $\gamma \in \Gamma'$. $\per(\Gamma',X_f)$ is a subgroup of $X_f$. If $f$ is invertible in $l^1(\Gamma)$ then whenever $\Gamma'$ has finite index in $\Gamma$, $\per(\Gamma',X_f)$ is finite. For a sequence $\{\Gamma_i\}_{i=1}^\infty$ of finite index subgroups of $\Gamma$, let 
  $$\gr\big(\{\Gamma_i\}_{i=1}^\infty,X_f\big):=\limsup_{n\to\infty} \frac{\log(|\per(\Gamma_i,X_f)|)}{[\Gamma:\Gamma_i]}$$
   be the growth rate of the number of periodic points in $X_f$ with respect to $\{\Gamma_i\}_{i=1}^\infty$. We will prove
   
   \begin{thm}\label{thm:periodic}
   If $f \in \Z\Gamma$ is invertible in $l^1(\Gamma)$ and $\Gamma$ is non-amenable then for any sequence $\{\Gamma_i\}_{i=1}^\infty$ of finite-index normal subgroups of $\Gamma$ with $\lim_{i\to\infty}\Gamma_i=\{e\}$,
   $$\gr\big(\{\Gamma_i\}_{i=1}^\infty,X_f\big) = h\big(\{\Gamma_i\}_{i=1}^\infty, X_f, \mu_f\big).$$
   \end{thm}   
      
     By [DS07, corollary 5.3 and theorem 6.1], $\gr\big(\{\Gamma_i\}_{i=1}^\infty,X_f\big) = \log \det (f)$. Hence theorem \ref{thm:main} follows from theorem \ref{thm:periodic}.
      
Along the way, we prove other results which might be of independent interest. In \S \ref{sec:stronglyergodic} we prove that if $f\in \Z\Gamma$ is invertible in $l^1(\Gamma)$, then the action $\Gamma \cc (X_f,\mu_f)$ is weakly equivalent to a Bernoulli shift action in the sense of [Ke09]. From this, we conclude that if $\Gamma $ is non-amenable then $\Gamma \cc (X_f,\mu_f)$ is strongly ergodic.

It seems possible to generalize theorem \ref{thm:main} to the $f$-invariant (defined in [Bo09a]) by making use of its characterization as a variant of sofic entropy given in [Bo09c].

   \subsection{Organization}
    The next section introduces notation and defines some important maps between relevant spaces. \S \ref{sec:upperbound}  shows that $h(\{\Gamma_i\}_{i=1}^\infty, X_f,\mu_f) \le \gr\big(\{\Gamma_i\}_{i=1}^\infty,X_f\big)$. \S \ref{sec:alternative} gives an alternative formulation of entropy. We prove in \S \ref{sec:stronglyergodic} that the action $\Gamma \cc (X_f,\mu_f)$ is strongly ergodic. These results are used in \S \ref{sec:lowerbound} to prove the oppositive inequality which completes the proof of theorem \ref{thm:periodic}. Results in [DS07] are then used to finish the proof of theorem \ref{thm:main}.
        
    {\bf Acknowledgments}.
    
    I'd like to thank Doug Lind for introducing me to this subject, Russ Lyons for helpful conversations and Hanfeng Li for pointing out errors in a previous version.

   
    \section{Preliminaries: a generating function}
    
    The purpose of this section is to introduce notation and define an explicit generating function $\phi:X_f \to A$ for $X_f$.
    
%
    
    Given $g\in l^\infty(\Gamma)$ and $h\in l^1(\Gamma)$ we define the convolutions $g\cdot h$ and $h\cdot g$ in $l^\infty(\Gamma)$ by
  \begin{eqnarray*}
  (g\cdot h)(\gamma) &:=& \sum_{\beta \in \Gamma} g(\beta) h(\beta^{-1}\gamma) =\sum_{\beta \in \Gamma} g(\gamma\beta^{-1}) h(\beta) \\
    (h\cdot g)(\gamma) &:=& \sum_{\beta \in \Gamma} h(\beta) g(\beta^{-1}\gamma) =\sum_{\beta \in \Gamma} h(\gamma\beta^{-1}) g(\beta).   \end{eqnarray*}
 
Let $\Z\Gamma\subset l^\infty(\Gamma)$ be the set of all functions $g:\Gamma \to \Z$ such that $g(\gamma)=0$ for all but finitely many $\gamma \in \Gamma$. Given sets $A$ and $B$, $A^B$ denotes the set of all functions from $B$ to $A$.  Let $\T=\R/\Z$. If $g\in \Z\Gamma$  and $h\in \T^\Gamma$ then define $g\cdot h$ and $h\cdot g \in \T^\Gamma$ as above.
          
  The {\em adjoint} of an element $g \in \C^\Gamma$ is defined by $g^*(\gamma) = \overline{g(\gamma^{-1})}$.

    Since $f \in \Z\Gamma$ is invertible in $l^1(\Gamma)$, there is an element $f^{-1}\in l^1(\Gamma)$ such that 
    $$f\cdot f^{-1}=f^{-1}\cdot f =1_e$$
    where $1_e \in l^1(\Gamma)$ is the identity element: $1_e(\gamma)=1$ if $\gamma=e$ and $1_e(\gamma)=0$ otherwise. Note that $1_e\cdot h = h \cdot 1_e = h$ for any $h\in l^\infty(\Gamma)$.  By definition, $X_f = \{ h \in \T^\Gamma:~ h\cdot f^* = 0\}$.

      Let $\pi:\R \to \R/\Z=\T$ be the quotient map. Let $\pi^\Gamma: \R^\Gamma \to \T^\Gamma$ be the quotient map defined by $\pi^\Gamma(x)(\gamma)=\pi(x(\gamma))$. Let $l^\infty(\Gamma,\Z)$ be the set of $g\in l^\infty(\Gamma)$ with $g(\gamma) \in \Z$ for all $\gamma\in\Gamma$. Define $l^\infty(\Gamma,\R)$ similarly. Observe that for any $g\in l^\infty(\Gamma,\R)$ and $h\in \Z\Gamma$, $\pi^\Gamma(g\cdot h) = \pi^\Gamma(g)\cdot h$. 
        
     Let $\hf$ be the adjoint of $f^{-1}$. Since $(g\cdot  h)^*=h^*\cdot g^*$ for any $g,h\in \C^\Gamma$ for which this is well-defined, it follows that $\hf \cdot f^* = (f\cdot f^{-1})^* = 1_e$. Define $\xi: l^\infty(\Gamma,\Z) \to \T^\Gamma$ by $\xi(h)= \pi^\Gamma(h\cdot \hf)$. Observe that 
     $$\xi(h)\cdot f^* = \pi^\Gamma(h\cdot \hf \cdot f^*) =  \pi^\Gamma(h) = \bz$$
     where $\bz \in \T^\Gamma$ is the zero element: $\bz(\gamma)=0+\Z$ for all $\gamma\in \Gamma$. This implies that $\xi$ maps $l^\infty(\Gamma,\Z)$ into $X_f$. Moreover, it a group homomorphism.
     
    Fix an irrational number $\kappa \in [0,1]$.  For $x \in \T^\Gamma$, let $L(x) \in l^\infty(\Gamma)$ be the element with values $L(x)(\gamma) \in [-1+\kappa,\kappa)$ for all $\gamma \in \Gamma$ defined by $L(x)(\gamma) = x(\gamma) \mod 1$. $L(x)$ is a `lift' of $x$. The reason for choosing $\kappa$ irrational is that $\mu_f(\{x \in X_f:~ x(e)=\kappa \mod \Z\})=0$ where $\mu_f$ is the Haar probability measure on $X_f$. This will be useful later.
    
    Define $P(x) \in l^\infty(\Gamma)$ by $P(x)=L(x)\cdot f^*$. If $x\in X_f$ then $x \cdot f^*=\bz$. Therefore, $L(x)\cdot f^* \in l^\infty(\Gamma,\Z)$. Thus $P$ maps $X_f$ into $l^\infty(\Gamma,\Z)$. Since the support of $f$ is finite, there exists an integer $M>0$ such that for all $x\in X_f$, $P(x) \in [-M,M]^\Gamma$. Also, 
     $$\xi(P(x)) = \pi^\Gamma(P(x)\cdot \hf) = \pi^\Gamma( L(x)\cdot f^* \cdot \hf) = \pi^\Gamma(L(x))= x.$$
    This proves that $\xi$ is surjective. We define $\phi:X_f \to [-M,M]\cap \Z=:I_M$ by $\phi(x) = P(x)(e)$. The calculation above implies that $\phi$ is generating. So theorem \ref{thm:K} implies that for any  sequence $\{\Gamma_i\}_{i=1}^\infty$ of finite-index normal subgroups of $\Gamma$ with $\lim_{i\to\infty}\Gamma_i=\{e\}$,
    $$h\big(\{\Gamma_i\}_{i=1}^\infty, X_f,\mu_f\big) = h\big(\{\Gamma_i\}_{i=1}^\infty, \phi\big).$$
    
   \subsection{Preliminaries: finite quotients}
   Let $\Gamma_i$ be a finite index subgroup of $\Gamma$.  It will be helpful to extend the definitions of $\xi,L,P$, etc. to functions of $\Gamma_i\backslash \Gamma$. So for $g \in l^1(\Gamma)$ and $\psi \in \R^{\Gamma_i\backslash \Gamma}$ define $\psi \cdot g \in \R^{\Gamma_i\backslash \Gamma}$  by
   $$\psi \cdot g(\Gamma_i \gamma):= \sum_{\beta \in \Gamma} \psi(\Gamma_i \gamma \beta^{-1})g(\beta).$$
   Similarly, if $g\in \Z\Gamma$ and $\psi \in \T^{\Gamma_i\backslash \Gamma}$ define $\psi \cdot g \in \T^{\Gamma_i\backslash \Gamma}$ as above. Let $\pi^{\Gamma_i\backslash \Gamma}:\R^{\Gamma_i\backslash \Gamma} \to \T^{\Gamma_i\backslash \Gamma}$ be the projection map $\pi^{\Gamma_i\backslash \Gamma}(\psi)(\Gamma_i \gamma) = \pi(\psi(\Gamma_i\gamma))$. Define $\xi:\Z^{\Gamma_i\backslash \Gamma} \to \T^{\Gamma_i \backslash \Gamma}$ by $\xi(\psi) = \pi^{\Gamma_i\backslash \Gamma}(\psi \cdot \hf)$. 
    

  If $\psi \in \T^{\Gamma_i\backslash \Gamma}$ then we define $L(\psi) \in \R^{\Gamma_i\backslash \Gamma}$ by $L(\psi)(\Gamma_i \gamma) \in [-1+\kappa,\kappa)$ and $L(\psi)(\Gamma_i \gamma) = \psi(\Gamma_i \gamma) \mod 1$ for all $\gamma \in \Gamma$. $P(\psi) \in \Z^{\Gamma_i\backslash \Gamma}$ is defined by $P(\psi):=L(\psi)\cdot f^*$.

   \section{Entropy and periodic points: upper bound }\label{sec:upperbound}
   
   The purpose of this section is to prove:
   
   \begin{thm}\label{thm:upperbound}
   If $\{\Gamma_i\}_{i=1}^\infty$ is a  sequence of finite-index normal subgroups of $\Gamma$ with $\lim_{i\to\infty}\Gamma_i=\{e\}$ then 
   $$\gr\big(\{\Gamma_i\}_{i=1}^\infty,X_f\big) \ge h\big(\{\Gamma_i\}_{i=1}^\infty,X_f,\mu_f\big).$$
    \end{thm}
       
    \begin{proof}
Here is a brief sketch of the proof.  Let $X_f^i:=\{x \in \T^{\Gamma_i\backslash\Gamma}:~x\cdot f^*=\bz\}$. Note that $\xi$ maps $\Z^{\Gamma_i\backslash \Gamma}$ to $X^i_f$. There is an obvious bijection between $X^i_f$ and $\per(\Gamma_i,X_f)$. So it suffices to bound, for fixed $x\in X^i_f$, the cardinality of the set
  $$\big\{\psi\in I_M^{\Gamma_i\backslash \Gamma}:~ \xi(\psi)=x, \textrm{ and } ||\phi^W_* \mu-\psi^W_*\zeta_i||_1<\epsilon\big\}$$
  where $W \subset \Gamma$ is a finite set of our choice and $I_M=[-M,M]\cap \Z$. To do this, we will show that if $\epsilon>0$ and $W$ is chosen appropriately, then $||\phi^W_* \mu-\psi^W_*\zeta_i||_1<\epsilon$ implies $P(\xi(\psi))(C)=\psi(C)$ for most $C \in   \Gamma_i\backslash \Gamma$.  
  
    Let $\epsilon>0$. We assume that $2\epsilon |\supp(f)| < 1$ where $\supp(f):=\{\gamma \in \Gamma:~f(\gamma)\ne 0\}$ denotes the support of $f$. Let $\rho:X_f \to \T$ be the evaluation map $\rho(x)=x(e)$. It is a continuous group homomorphism. So $\rho_*\mu_f$ is the Haar measure of a closed subgroup of $\T$. Since $\kappa$ is irrational, there exists a  $\delta>0$ such that 
    $$\rho_*\mu_f\big(\{t\in \T:~ t \in [\kappa-2\delta,\kappa+2\delta] \mod 1\}\big)< \epsilon.$$
  
   
   For any integer $n>0$, let $I_n = [-n,n] \cap \Z$. Recall that $M>0$ is such that $P(X_f) \subset (I_M)^\Gamma=:I_M^\Gamma$. Because $\hf =(f^{-1})^*\in l^1(\Gamma)$, there exists a finite set $W \subset \Gamma$ such that 
    $$M\sum_{\gamma \in \Gamma - W} |\hf(\gamma)| < \delta.$$
    Let 
    $$\cB:=\Big\{g\in (I_M)^{W}:~ \sum_{w \in W} g(w) \hf(w) \notin (-1+\kappa+\delta,\kappa-\delta) \Big\}.$$
{\bf Claim 1}. $\phi^{W}_*\mu_f(\cB)<\epsilon$. 
\begin{proof}
If $x \in X_f$ and $\phi^{W}(x) \in \cB$ then $\sum_{w\in W} \phi^{W}(x)(w)\hf(w) \notin  (-1+\kappa+\delta,\kappa-\delta)$. But 
   $$\phi^{W}(x)(w) = \phi(wx) = P(wx)(e)=P(x)(w^{-1}).$$
   So,  $\sum_{w\in W} P(x)(w^{-1})\hf(w) \notin  (-1+\kappa+\delta,\kappa-\delta)$     
   By definition of $W$, this implies 
   $$P(x)\cdot \hf(e)=\sum_{\gamma \in \Gamma} P(x)(\gamma^{-1})\hf(\gamma) \notin  (-1+\kappa+2\delta,\kappa-2\delta).$$
  Since $P(x)\cdot \hf = L(x) \in [-1+\kappa,\kappa)^\Gamma$, it follows that $x(e) \notin  (-1+\kappa+2\delta,\kappa-2\delta) \mod 1$. The choice of $\delta$ now implies the claim.
\end{proof}

Let  $\zeta_i$ be the uniform probability measure on $\Gamma_i \backslash \Gamma$.  Claim 1 implies that if $\psi \in I_M^{\Gamma_i\backslash \Gamma}$  and $||\phi^W_* \mu_f-\psi^W_*\zeta_i||_1<\epsilon$, then $\psi^W_*\zeta_i(\cB) \le 2\epsilon$. 
   
   For $C\in \Gamma_i\backslash \Gamma$ and $\psi\in I_M^{\Gamma_i\backslash \Gamma}$, define $\psi_C \in I_M^W$ by $\psi_C(w):=\psi(Cw^{-1})$. Let $\cV(\psi):=\{C\in \Gamma_i\backslash \Gamma:~ \psi_C \in \cB\}.$ The previous paragraph implies that if  $||\phi^W_* \mu_f-\psi^W_*\zeta_i||_1<\epsilon$ then $|\cV(\psi)| \le 2\epsilon[\Gamma:\Gamma_i]$.   

\noindent {\bf Claim 2}. If $C\notin \cV(\psi)$ then $L(\xi(\psi))(C)=\psi\cdot \hf(C)$. 

\begin{proof}
By definition of $\cV$, $\sum_{w \in W} \psi_C(w) \hf(w)  \in (-1+\kappa+\delta,\kappa-\delta)$. Note that:
  $$\psi\cdot \hf(C) = \sum_{\gamma \in \Gamma} \psi(C\gamma^{-1})\hf(\gamma) = \sum_{w \in W} \psi_C(w) \hf(w)  +   \sum_{\gamma \in \Gamma-W} \psi(C\gamma^{-1})\hf(\gamma).$$ 
    The definition of $W$ now implies $\psi\cdot \hf(C) \in (-1+\kappa,\kappa)$. The definitions of $L$ and $\xi$ imply $L(\xi(\psi))(C)=\psi\cdot\hf(C)$ as claimed. 
\end{proof}
   
  Define
   $$\cV'(\psi) :=\bigcup_{\gamma \in \supp(f)} \cV(\psi)\gamma^{-1} = \big\{C \in \Gamma_i\backslash \Gamma:~\exists \gamma \in \supp(f) \textrm{ such that } C\gamma \in \cV(\psi) \big\}.$$
   If $||\phi^W_* \mu_f-\psi^W_*\zeta_i||_1<\epsilon$ then since $|\cV(\psi)| \le 2\epsilon[\Gamma:\Gamma_i]$, we must have $|\cV'(\psi)| \le 2\epsilon |\supp(f)| [\Gamma:\Gamma_i]$. Claim 2 implies that if $C\notin \cV'(\psi)$ then $L(\xi(\psi))(C\gamma) = \psi\cdot\hf(C\gamma)$ for all $\gamma \in \supp(f)$. Thus,
   \begin{eqnarray*}
   P(\xi(\psi))(C) &=& L(\xi(\psi))\cdot f^*(C)
                        = \sum_{\gamma \in \Gamma} L(\xi(\psi))(C\gamma)f^*(\gamma^{-1})\\
                        &=& \sum_{\gamma \in \supp(f)} L(\xi(\psi))(C\gamma)f(\gamma)
                        = \sum_{\gamma \in \supp(f)} \psi\cdot\hf(C\gamma)f^*(\gamma^{-1})\\
                        &=& \psi\cdot \hf \cdot f^*(C) = \psi(C).
   \end{eqnarray*}
   
   We have proven above that  if $x \in X^i_f$, $\psi \in I_M^{\Gamma_i\backslash\Gamma}$, $\xi(\psi)=x$ and $||\phi^W_* \mu_f-\psi^W_*\zeta_i||_1<\epsilon$ then there exists a set $\cV'(\psi)$ of cardinality $|\cV'(\psi)| \le 2\epsilon |\supp(f)| [\Gamma:\Gamma_i]$ such that $\psi(C)=P(x)(C)$ for all $C\notin \cV'(\psi)$. Thus for any $x \in X^i_f$,
   $$\#\big\{\psi \in    I_M^{\Gamma_i\backslash\Gamma}: ~      ||\phi^W_* \mu_f-\psi^W_*\zeta_i||_1<\epsilon, ~\xi(\psi)=x\big\}\le {[\Gamma:\Gamma_i] \choose  \lfloor 2\epsilon |\supp(f)| [\Gamma:\Gamma_i]\rfloor} (2M+1)^{2 \epsilon |\supp(f)| [\Gamma:\Gamma_i]}.$$
   Since $|X^i_f|=|\per(\Gamma_i,X_f)|$, this implies
$$\#\big\{\psi \in    I_M^{\Gamma_i\backslash\Gamma}: ~     ||\phi^W_* \mu_f-\psi^W_*\zeta_i||_1<\epsilon\big\} \le   |\per(\Gamma_i,X_f)|   {[\Gamma:\Gamma_i] \choose  \lfloor2\epsilon |\supp(f)| [\Gamma:\Gamma_i]\rfloor} (2M+1)^{2\epsilon |\supp(f)| [\Gamma:\Gamma_i]}.$$
Stirling's approximation now implies
    $$ h\big(\{\Gamma_i\}_{i=1}^\infty,X_f,\mu_f\big) \le \gr\big(\{\Gamma_i\}_{i=1}^\infty,X_f\big) + H\big( 2\epsilon|\supp(f)| \big) + 2\epsilon|\supp(f)|\log(2M+1)$$
    where $H(x):=-x\log(x)-(1-x)\log(1-x)$. Since $\epsilon>0$ is arbitrary this proves     $ h\big(\{\Gamma_i\}_{i=1}^\infty,X_f,\mu_f\big) \le \gr\big(\{\Gamma_i\}_{i=1}^\infty,X_f\big)$ and completes the theorem.    
     \end{proof}
  
  \begin{remark}
  The above proof does not use the fact that $\{\Gamma_i\}_{i=1}^\infty$ is a  sequence of finite-index normal subgroups of $\Gamma$ with $\lim_{i\to\infty}\Gamma_i=\{e\}$. Indeed, the theorem holds true as long as each $\Gamma_i$ has finite-index in $\Gamma$ and $\lim_{i\to\infty} [\Gamma:\Gamma_i] = +\infty$.
  \end{remark}

  \section{An alternative formula for entropy}\label{sec:alternative}
  
  The purpose of this section is to provide an alternative formula for the entropy of an ergodic action. This will be used in \S \ref{sec:lowerbound} to prove a lower bound on entropy in terms of the growth rate of periodic points.
  
  \begin{defn}
  Given a finite set $A$, let $A^\Gamma$ denote the set of all functions $x:\Gamma \to A$. We let $A^\Gamma$ have the topology of uniform convergence on finite subsets. $\Gamma$ acts on $A^\Gamma$ by
   $$(\gamma x)(g)=x(\gamma^{-1}g) ~\forall \gamma,g \in \Gamma, x\in A^\Gamma.$$
  We will say that a Borel measure $\mu$ on $A^\Gamma$ is {\em shift-invariant} if for 
  $$\mu(\gamma E) = \mu(E) ~\forall \textrm{ Borel } E \subset A^\Gamma, \gamma \in \Gamma.$$
  Let  $M(A^\Gamma)$ be the space of all shift-invariant Borel probability measures on $A^\Gamma$. The {\em weak* topology} on $M(A^\Gamma)$ is defined as follows. A sequence $\{\mu_i\}_{i=1}^\infty \subset M(A^\Gamma)$ converges to a measure $\mu_\infty$ if and only if $\lim_{i\to\infty} \int c ~d\mu_i  = \int c ~d\mu_\infty$ for every continuous function $c$ on $A^\Gamma$. By the Banach-Alaoglu theorem, $M(A^\Gamma)$ is sequentially compact in the weak* topology.      \end{defn}

 \begin{defn}
  Let $\Gamma$ act on a standard probability space $(X,\mu)$ by measure-preserving transformations. Let $\phi:X \to A$ be a measurable map into a finite set. Then $\phi$ induces a map $\tphi$ from $X$ to $A^\Gamma$ by $\tphi(x)(\gamma)=\phi(\gamma^{-1} x)$. This map is equivariant. Therefore, $\tphi_*\mu$ is a shift-invariant probability measure on $A^\Gamma$.
  \end{defn}
  
  
  \begin{defn}
  If $\Gamma'$ is a subgroup of $\Gamma$, then let $M(\Gamma', A^\Gamma)$ be the set of all measures $\eta \in M(A^\Gamma)$ such that $\eta$ is supported on $\per(\Gamma',A^\Gamma)$ which equals the set of all $x \in A^\Gamma$ such that $\gamma' x = x$ for all $\gamma'\in\Gamma'$.
  \end{defn}
  
  \begin{defn}
   Let $\{\Gamma_i\}_{i=1}^\infty$ be a  sequence of finite-index normal subgroups of $\Gamma$ with $\lim_{i\to\infty}\Gamma_i=\{e\}$.  Let $\{\mu_i\}_{i=1}^\infty$ be a sequence with $\mu_i \in M(\Gamma_i, A^\Gamma)$. Define $H(\mu_i)$ as usual:
   $$H(\mu_i) := -\sum_{x \in A^\Gamma} \mu_i\big(\{x\}\big)\log\big(\mu_i(\{x\})\big).$$
   Let $$h\big( \{\mu_i\}_{i=1}^\infty\big):=\limsup_i \frac{H(\mu_i)}{[\Gamma:\Gamma_i]}.$$
   Let $\bh\big(\{\Gamma_i\}_{i=1}^\infty,\phi\big):= \sup h\big( \{\mu_i\}_{i=1}^\infty\big)$ where the $\sup$ is over all sequences $\{\mu_i\}_{i=1}^\infty$ with $\mu_i \in M(\Gamma_i, A^\Gamma)$ such that $\{\mu_i\}_{i=1}^\infty$ converges to $\tphi_*\mu$ in the weak* topology.
      \end{defn}
 
 The purpose of this section is to prove the following.
 \begin{thm}\label{thm:alternative}
 If $\Gamma \cc (X,\mu)$ is ergodic and $\{\Gamma_i\}_{i=1}^\infty$ is a  sequence of finite-index normal subgroups of $\Gamma$ with $\lim_{i\to\infty}\Gamma_i=\{e\}$, then for any measurable map $\phi:X \to A$ into a finite set $A$,  $\bh\big(\{\Gamma_i\}_{i=1}^\infty,\phi\big)=h\big(\{\Gamma_i\}_{i=1}^\infty,\phi\big)$.
  \end{thm}

  \begin{lem}
  With notation as above, $\bh\big(\{\Gamma_i\}_{i=1}^\infty,\phi\big)\ge h\big(\{\Gamma_i\}_{i=1}^\infty,\phi\big)$.
    \end{lem}
  \begin{proof}
  Without loss of generality, we assume that $h(\{\Gamma_i\}_{i=1}^\infty,\phi) \ne -\infty$.  For $W \subset G$ finite, $\epsilon>0$ and $i\in \N$, let $\nu_{W,\epsilon,i}$ be the uniform probability measure on the set $\{\psi\in \Gamma_i\backslash \Gamma \to A:~ ||\phi^W_* \mu-\psi^W_*\zeta_i||_1<\epsilon\}$ where $\zeta_i$ is the uniform probability measure on $\Gamma_i\backslash \Gamma$. The assumption $h(\{\Gamma_i\}_{i=1}^\infty,\phi) \ne -\infty$ implies this set is nonempty if $i$ is sufficiently large.
  
  Let $\Omega:\per(\Gamma_i,A^\Gamma) \to A^{\Gamma_i\backslash\Gamma}$ be the function $\Omega(x)(\Gamma_i g) := x(g)$. $\Omega$ is a bijection. Let $\tnu_{W,\epsilon,i}$ be the measure on $\per(\Gamma_i,A^\Gamma)$ defined by $\tnu_{W,\epsilon,i}(E)=\nu_{W,\epsilon,i}(\Omega(E))$. Then 
$$\log\Big( \#\big\{\psi:\Gamma_i \backslash \Gamma \to A:~||\phi^W_* \mu-\psi^W_*\zeta_i||_1<\epsilon\big\}\Big) = H(\tnu_{W,\epsilon,i}).$$
Observe that the quantity above decreases as $\epsilon\to 0$. Also, it is monotone decreasing in $W$ where the set of finite subsets of $\Gamma$ is ordered by inclusion. So there exist sequences $\{W_i\}_{i=1}^\infty$ and $\{\epsilon_i\}_{i=1}^\infty$ such that $W_i$ is increasing, $\cup_i W_i =\Gamma$ and $\epsilon_i \searrow 0$ and if $\tnu_i:=\tnu_{W_i,\epsilon_i,i}$ then $h(\{\tnu_i \}_{i=1}^\infty) = h\big(\{\Gamma_i\}_{i=1}^\infty,\phi\big)$. 

We claim that $\big\{\tnu_i\big\}_{i=1}^\infty$ converges to $\tphi_*\mu$. To see this note that if $E:A^\Gamma \to A$ is the evaluation map $E(x):=x(e)$ then for any finite $W\subset \Gamma$, $E^W_*(\tphi_*\mu)=\phi^W_*\mu$. Also, if $i$ is sufficiently large then $W_i \supset W$. Therefore, $E^W_*\tnu_i$ is a convex sum of measures of the form $\psi^W_*\zeta_i$ where $\psi:\Gamma_i\backslash \Gamma\to A$ is such that  $||\phi^{W_i}_* \mu-\psi^{W_i}_*\zeta_i||_1<\epsilon_i$. Thus we must have $||E^W_*(\tphi_*\mu) - E^W_*\nu_i||_1<\epsilon_i$. Since $\epsilon_i \to 0$ as $i\to\infty$ and $W \subset\Gamma$ is arbitrary, this implies the claim.

We now have
$$\bh\big(\{\Gamma_i\}_{i=1}^\infty,\phi\big)\ge h\big(\{\tnu_i\}_{i=1}^\infty\big) = h\big(\{\Gamma_i\}_{i=1}^\infty,\phi\big)$$
as claimed.
  \end{proof}
  
  \begin{proof}[Proof of theorem \ref{thm:alternative}]
  
   Let $\{\eta_i\}_{i=1}^\infty$ be a sequence of probability measures with $\eta_i \in M(A^\Gamma,\Gamma_i)$. Suppose $\{\eta_i\}_{i=1}^\infty $ converges to $\tphi_*\mu$ in the weak* topology. Let $ W \subset G$ be finite and $\epsilon>0$.  By the previous lemma, it suffices to show that 
$$h\big(\{\eta_i\}\big) \le \limsup_{i\to\infty}  \frac{\log\big( \#\{\psi\in A^{\Gamma_i \backslash \Gamma}:~||\phi^W_* \mu-\psi^W_*\zeta_i||_1<\epsilon\}\big)}{[\Gamma:\Gamma_i]}.$$
Let $\zeta_i$ be the uniform probability measure on $\Gamma_i \backslash \Gamma$. For every function $s:A^W \to \{-1,+1\}$, let 
$$B_i(s):=\Big\{ \psi\in A^{\Gamma_i \backslash \Gamma} :~ ||\phi^W_* \mu-\psi^W_*\zeta_i||_1 \ge \epsilon \textrm{ and }  \forall j\in A^W,~s(j)\big(\psi^W_*\zeta_i(\{j\}) -\phi^W_*\mu(\{j\})\big) \ge 0\Big\}.$$ 
$s$ is for {\em sign} and $B$ is for {\em bad}. Let
$$G_i:=\big\{\psi \in A^{\Gamma_i \backslash \Gamma}:~||\phi^W_* \mu-\psi^W_*\zeta_i||_1 < \epsilon\big\}.$$
Observe that $G_i$ and $\{B_i(s):~s\in \{-1,+1\}^{A^W} \}$ partition $A^{\Gamma_i \backslash \Gamma}$. As in the previous lemma, let $\Omega:\per(\Gamma_i,A^\Gamma) \to A^{\Gamma_i\backslash\Gamma}$ be the function $\Omega(x)(\Gamma_i g) := x(g)$. $\Omega$ is a bijection. Let $\tB_i(s) = \Omega^{-1}(B_i(s))$ and $\tG_i = \Omega^{-1}(G_i)$. Because each $\Gamma_i<\Gamma$ is normal, these sets are $\Gamma$-invariant: $\gamma \tB_i(s)=\tB_i(s)$ and $\gamma \tG_i=\tG_i$ for all $\gamma \in \Gamma$. 

 If $\eta_i(\tB_i(s))>0$ then let $\eta_{i,s}$ be the measure on $A^{\Gamma_i \backslash \Gamma}$ defined by
$$\eta_{i,s}(E):=\frac{\eta_i(E \cap \tB_i(s))}{\eta_i(\tB_i(s))}.$$
Otherwise, choose $\eta_{i,s} \in M(A^\Gamma)$ arbitrarily. Similarly, if $\eta_i(\tG_i)>0$ then let $\eta_{i,G}:=\frac{\eta_i(E\cap \tG_i)}{\eta_i(\tG_i)}$ and otherwise choose $\eta_{i,G} \in M(A^\Gamma)$ arbitrarily. Thus if $t_{i,s} = \eta_i(\tB_i(s))$ and $t_{i,G}=\eta_i(\tG_i)$ then
$$\eta_i := t_{i,G}\eta_{i,G} + \sum_{s:A^W \to \{-1,+1\}} t_{i,s} \eta_{i,s}.$$
Because $M(A^\Gamma)$ is sequentially compact, after passing to a subsequence if necessary, we may assume that $\{\eta_{i,G}\}_{i=1}^\infty$ converges in the weak* topology to a measure $\omega_G$ and for each $s:A^W \to \{-1,+1\}$, $\{\eta_{i,s}\}_{i=1}^\infty$ converges in the weak* topology to a measure $\omega_s$. We may also assume that there are real numbers $r_G, r_s$ (for $s: A^W \to \{-1,+1\})$ such that $\lim_i t_{i,G} = r_G$ and $\lim_i t_{i,s} = r_s$. Since $\{\eta_i\}_{i=1}^\infty$ converges to $\tphi_*\mu$, 
$$\tphi_*\mu = r_G \omega_G + \sum_{s:A^W \to \{-1,+1\}} r_{s} \omega_{s}.$$  
By construction, for any $s: A^W \to \{-1,+1\}$ if $r_s \ne 0$ then $\omega_s$ cannot equal $\tphi_*\mu$. Indeed if $E: A^\Gamma \to A$ is the map $E(x)=x(e)$ then $E^W_*\eta_{i,s}$ converges to $E^W_*\omega_s$. By construction, for any $j\in A^W$,
$$s(j)\Big( E^W_*\eta_{i,s}(\{j\}) - E^W_*\tphi_*\mu(\{j\})\Big) \ge 0$$
and
$$\sum_{j \in A^W} \Big| E^W_*\eta_{i,s}(\{j\})  - E^W_*\tphi_*\mu(\{j\})\Big|\ge \epsilon.$$
Hence,
$$|| E^W_*\omega_s - E^W_*\tphi_*\mu||_1 = \sum_{j \in A^W}  s(j)\Big( E^W_*\omega_s(\{j\})  - E^W_*\tphi_*\mu(\{j\})\Big)\ge \epsilon.$$
Since $\tphi_*\mu$ is ergodic, this implies that $r_s$ equals $0$ for each $s$. Thus we have proven that for all $s$, $\lim_i \eta_i(\tB(s)) = 0$. 

Since $\eta_i := t_{i,G}\eta_{i,G} + \sum_{s:A^W \to \{-1,+1\}} t_{i,s} \eta_{i,s}$ and the supports of $\eta_{i,G}$ and $\eta_{i,s}$ are disjoint,
$$H(\eta_i) = t_{i,G}\Big( H(\eta_{i,G}) - \log(t_{i,G})\Big) + \sum_{s:A^W \to \{-1,+1\}} t_{i,s} \Big(H(\eta_{i,s})- \log(t_{i,s})\Big).$$
Because $t_{i,G} \to 1$ and $t_{i,s} \to 0$ as $i \to \infty$ and $\frac{H(\eta_{i,s})}{[\Gamma:\Gamma_i]} \le \log(2M+1)$, it follows that 
$$h\big(\{\eta_i\}\big)=\limsup_i \frac{H(\eta_i)}{[\Gamma:\Gamma_i]} = \limsup_i \frac{H(\eta_{i,G})}{[\Gamma:\Gamma_i]}=h\big(\{\eta_{i,G}\}\big).$$
  However, $$H(\eta_{i,G}) \le \log(|G_i|) = \log\big( \#\{\psi: \Gamma_i \backslash \Gamma\to A:~||\phi^W_* \mu-\psi^W_*\zeta_i||_1<\epsilon\}\big).$$
So this proves the theorem.
  \end{proof}

  \section{Strong ergodicity}\label{sec:stronglyergodic}
 Here is a brief sketch of this section. In order to apply the results of the previous section, we need to show that the action $\Gamma \cc (X_f,\mu_f)$ is ergodic. Rather than solve this problem in general, we prove that if $\Gamma$ is non-amenable then this action is strongly ergodic, and therefore is ergodic. To do this, we represent $\mu_f$ as a weak* limit of Bernoulli factors.  From this characterization, we show that $\Gamma \cc (X_f,\mu_f)$ is weakly equivalent to a Bernoulli shift action in the sense of [Ke09]. Since Bernoulli shift actions are known to be strongly ergodic (whenever $\Gamma$ is non-amenable) this shows that $\Gamma \cc (X_f,\mu_f)$ is strongly ergodic. 
  

\begin{lem}\label{lem:Bernoulli}
For $n \in \Z$, $n>0$, let $u_n$ be the uniform measure on $I_n:=[-n,n]\cap\Z$. Let $\nu_n$ be the product measure $\nu_n:=u_n^\Gamma$ on $I_n^\Gamma$. Then $\{\xi_*\nu_n\}$ converges to $\mu_f$ as $n \to \infty$.
\end{lem}                         
 \begin{proof}
Let $W\subset \Gamma$ be finite. Let $x \in I_M^\Gamma$ have support in $W$. So, $x(\gamma)=0$ if $\gamma \notin W$.   Let $\omega$ be any weak* limit point of  $\{\xi_*\nu_n\}_{n=1}^\infty$. We claim that $\omega$ is invariant under addition by $\xi(x)$. That is, if $A \subset X_f$ is Borel then $\omega(A) = \omega(A + \xi(x))$. To see this, let
$$Z_n:=\{ y\in I_n^\Gamma:~ y(w) \in I_{n-M} ~\forall w\in W\}.$$
Note that if $B\subset Z_n$ is Borel then $\nu_n(B) = \nu_n(B+x)$. Also note that $\nu_n(Z_n)=\Big(\frac{2n-2M+1}{2n+1}\Big)^{|W|}$ tends to $1$ as $n \to \infty$. 


Now let $c:X_f \to \R$ be a continuous function. 
\begin{eqnarray*}
\Big|\int c\big(y+\xi(x)\big) -c(y)~d\xi_*\nu_n(y) \Big| &=&\Big|\int c\big(\xi(z)+\xi(x)\big) -c\big(\xi(z)\big)~d\nu_n(z) \Big| \\
&\le&\Big|\int_{Z_n} c\big(\xi(z)+\xi(x)\big) -c\big(\xi(z)\big)~d\nu_n(z) \Big| + 2\big(1-\nu_n(Z_n)\big)||c||_\infty\\
&=& 2\big(1-\nu_n(Z_n)\big)||c||_\infty.
 \end{eqnarray*}
 Since $\nu_n(Z_n)$ tends to $1$ as $n\to\infty$ and $\omega$ is a weak* limit point of $\{\xi_*\nu_n\}_{n=1}^\infty$, 
 $$\int c(y+\xi(x)) ~d\omega(y) =\int c(y) ~d\omega(y).$$
Since $c$ is arbitrary, this implies that $\omega$ is invariant under addition by $\xi(x)$ (i.e., $\omega(E+\xi(x))=\omega(E)$ for any Borel $E \subset X_f$). The set $\{\xi(x):~ x \in I_M^\Gamma$ has finite support $\}$ is dense in $X_f$ (since every element of $X_f$ equals $\xi(y)$ for some $y\in I_M^\Gamma$ and $\xi$ is continuous). Therefore, $\omega$ is invariant under addition by all elements of $X_f$. By uniqueness of Haar measure, $\omega=\mu_f$.
\end{proof}


\begin{defn}
Let $\Gamma \cc (Y,\nu)$ and $\Gamma \cc (Z,\zeta)$ be two measure-preserving actions of $\Gamma$ on standard probability spaces. We say that  $\Gamma \cc (Y,\nu)$ {\em weakly contains} $\Gamma \cc (Z,\zeta)$ if for every measurable map $\psi:Z \to A$ (where $A$ is a finite set), and every finite $W\subset \Gamma$ and every $\epsilon>0$ there exists a measurable map $\phi:Y \to A$ such that 
$$ \big|\big| \phi^W_*\nu - \psi^W_*\zeta \big|\big|_1 \le \epsilon.$$

If $\Gamma \cc (Y,\nu)$  weakly contains $\Gamma \cc (Z,\zeta)$ and  $\Gamma \cc (Z,\zeta)$  weakly contains $\Gamma \cc (Y,\nu)$ then we say $\Gamma \cc (Y,\nu)$ is {\em weakly equivalent} to $\Gamma \cc (Z,\zeta)$. These notions were introduced in [Ke09] (in slightly different language).
\end{defn}
We will show that $\Gamma \cc (X_f,\mu_f)$ is weakly equivalent to a Bernoulli shift action. To do this, we need to understand how weak containment behaves under weak* limits. This is handled in the next proposition.

\begin{prop}\label{prop:weak*}
Let $X$ be a compact metric space. Suppose $\Gamma$ acts on $X$ by homeomorphisms. Let $\Gamma \cc (Y,\nu)$ be a measure-preserving action on a standard probability space. Then the set of all $\Gamma$-invariant Borel probability measures $\omega$ on $X$ such that $\Gamma \cc (Y,\nu)$ weakly contains $\Gamma \cc (X,\omega)$ is closed in the weak* topology.
\end{prop}

\begin{proof}
Let $\{\mu_i\}_{i=1}^\infty$ be a sequence of $\Gamma$-invariant Borel probability measures on $X$ that converge in the weak* topology to a measure $\mu_\infty$. Suppose that for each $i\ge 1$, $\Gamma \cc (X,\mu_i)$ is weakly contained in $\Gamma \cc (Y,\nu)$. It suffices to show that $\Gamma \cc (X,\mu_\infty)$ is weakly contained in $\Gamma \cc (Y,\nu)$.

Let $\phi:X \to A$ be a $\mu_\infty$-measurable map into a finite set $A$. Because $\mu_\infty$ is a Borel measure, we may assume, after changing $\phi$ on a set of $\mu_\infty$-measure zero if necessary, that $\phi$ is Borel. Let $\epsilon>0$ and $W\subset \Gamma$ be finite. Choose $\delta>0$ so that  $|A^W|\big(|W|\delta + 3|A|^2|W|\delta^{1/4}\big)  +\delta < \epsilon$ and $|A|^{-1}>\delta^{1/4}$. This will be useful later.

For each $a\in A$ let $E_a=\phi^{-1}(a)$.  Because $X$ is a compact metric space, there exists a closed set $F_a \subset E_a$ and an open set $O_a \supset E_a$ such that $\mu_\infty(O_a-F_a)<\delta$. After replacing $O_a$ with $O_a \cap \bigcap_{b \ne a} (X-F_b)$ if necessary, we may assume that  $O_a \cap F_b = \emptyset$ if $a\ne b$. A standard partitions of unity argument (see e.g., [Ro88, proposition 9.16]) implies that there exists a collection of continuous function $\{C_a\}_{a\in A}$ such that 
\begin{itemize}
\item $C_a:X \to [0,1]$,
\item $\{x \in X:~ C_a(x)>0\} \subset O_a$,
\item $\forall x\in X$, $\sum_{a\in A} C_a(x) = 1$.
\end{itemize}

For $a \in A$,  let $\chi_a:X \to \R$ be the characteristic function defined by $\chi_a(x)=1$ if $\phi(x)=a$ and $\chi_a(x)=0$ otherwise.

\noindent {\bf Claim 1}. For any $a\in A$, $\int |C_a(x) - \chi_a(x)|~d\mu_\infty < \delta$.

\begin{proof}  Since $\sum_a C_a(x)=1$ for any $x\in X$ and $O_a \cap F_b = \emptyset$ if $a\ne b$, we must have $C_a(x)=1$ for all $x\in F_a$. Thus $C_a(x)=\chi_a(x)$ for $x\in F_a$. Since both $C_a$ and $\chi_a$ vanish outside of $O_a$, this implies
$$\int |C_a(x) - \chi_a(x)|~d\mu_\infty = \int_{O(a)-F(a)} |C_a(x) - \chi_a(x)|~d\mu_\infty <\delta.$$
 \end{proof}

\noindent {\bf Claim 2}. If $a,b \in A$ and $a\ne b$ then $\int C_aC_b ~d\mu_\infty < \delta$.

\begin{proof} By definition, $C_a$ vanishes on $X-O_a$ and $C_b$ vanishes on $X-O_b \supset F_a$. Hence
$$\int C_aC_b ~d\mu_\infty = \int_{O_a-F_a} C_aC_b ~d\mu_\infty <\delta.$$
\end{proof}
For any $j\in A^W$, define
$$C_j(x):=\prod_{w \in W} C_{j(w)} \big( wx\big).$$
Since $\{\mu_i\}_{i=1}^\infty$ converges to $\mu_\infty$ in the weak* topology, there exists an $N>0$ such that $n>N$ implies 
\begin{eqnarray*}
\Big|\int C_aC_b ~d\mu_n \Big|&<&\delta ~\forall a\ne b \in A\\
\sum_{j\in A^W} \Big| \int C_j ~d\mu_\infty - \int C_j ~d\mu_n\Big| &<& \delta.
\end{eqnarray*}

Observe that for any $a\ne b \in A$,
\begin{eqnarray}\label{eqn:thing1}
\mu_n\big( \{x:~ C_a(x)\ge \delta^{1/4} \textrm{ and } C_b(x) \ge \delta^{1/4}\} \big) &\le& \frac{1}{\sqrt{\delta}} \int C_aC_b ~d\mu_n< \sqrt{\delta}.
\end{eqnarray}

For each $a \in A$, let $\chi'_a:X\to \R$ be the function
\begin{displaymath}
\chi'_a(x)=\left\{\begin{array}{ll}
 1 &\textrm{ if } C_a(x) \ge \delta^{1/4} \textrm{ and } C_b(x) < \delta^{1/4}~\forall b \ne a\\
 0 & \textrm{ otherwise}
 \end{array}\right.
 \end{displaymath}
Since $|A|^{-1}>\delta^{1/4}$, it follows that for any $x\in X$, there exists some $a\in A$ such that $C_a(x)>\delta^{1/4}$. So equation (\ref{eqn:thing1}) implies: 
\begin{eqnarray*}
\mu_n\big( \{x:~ \chi'_a(x)=0 ~\forall a\in A\}) &\le& \sum_{a\ne b \in A}\mu_n\big( \{x:~ C_a(x)\ge \delta^{1/4} \textrm{ and } C_b(x) \ge \delta^{1/4}\} \big)\le |A|^2\sqrt{\delta}.
\end{eqnarray*}

Let $\psi:X \to A$ be any Borel function such that $\psi(x)=a$ if the maximum of $\{C_b(x):~b \in A\}$ is uniquely attained by $C_a(x)$. For each $a\in A$, let $\chi''_a$ be the characteristic function of $\psi^{-1}(a)$. So $\chi''_a \ge \chi'_a$. 

\noindent {\bf Claim 3}. For any $a\in A$, $\int |\chi''_a - C_a| ~d\mu_n < |A|^2\sqrt{\delta} +2|A|\delta^{1/4} \le 3|A|^2\delta^{1/4}$.
\begin{proof}
By definition, 
\begin{eqnarray}
\int |\chi''_a - C_a| ~d\mu_n  &\le&  \mu\big(\{x:~ \chi'_a(x)=0 ~\forall a\in A\}\big)+ \sum_{b\in A} \int \chi'_b |\chi''_a - C_a| ~d\mu_n\\
&\le& |A|^2\sqrt{\delta} + \sum_{b\in A} \int \chi'_b |\chi''_a - C_a| ~d\mu_n.\label{thing2}
\end{eqnarray}
If $b\ne a$ then $\chi'_b|\chi''_a-C_a|\le \delta^{1/4}$. Because $\sum_{a\in A} C_a(x)=1$ for any $x\in A$, it follows that if $\chi'_a(x)=1$ then $C_a(x) =1-\sum_{b\ne a}C_b(x) \ge 1-|A|\delta^{1/4}$. Hence,
$$\chi'_a|\chi''_a - C_a| \le |A|\delta^{1/4}.$$
This and equation (\ref{thing2}) proves the claim.
\end{proof}

For $j \in A^W$, define the functions $\chi''_j,\chi_j:X \to \R$ by
\begin{eqnarray*}
\chi''_j(x)&:=&\prod_{w \in W} \chi''_{j(w)}\big( wx\big)\\
\chi_j(x)&:=&\prod_{w \in W} \chi_{j(w)}\big( wx\big)\\
\end{eqnarray*}
Observe that $\chi_j$ is the characteristic function of $(\phi^W)^{-1}(j)$.

\noindent {\bf Claim 4}. For any $j \in A^W$,
\begin{eqnarray*}
\int |\chi_j - C_j| ~d\mu_\infty &\le & |W|\delta.
\end{eqnarray*}
\begin{proof}
We prove this by induction on $|W|$. If $|W|=1$ then the claim is true by claim 1. So assume $|W|>1$. Let $w \in W$ and $W' = W-\{w\}$. Let $j\in A^W$ and define $j' :W' \to A$ by $j'(w)=j(w) ~\forall w\in W'$. For any $x \in X$,
$$\chi_j(x) = \chi_{j'}(x) \chi_{j(w)}(wx) \textrm{ and }C_j(x) = C_{j'}(x) C_{j(w)}(wx).$$
Hence
$$|\chi_j(x) - C_j(x)| \le \Big|\big(\chi_{j'}(x)-C_{j'}(x)\big)\chi_{j(wx)}\Big| + \Big|C_{j'}(x)\big(\chi_j(wx) - C_{j(w)}(wx)\big)\Big|.$$
The induction hypothesis implies $\int |(\chi_{j'}(x)-C_{j'}(x))\chi_{j(wx)}| ~d\mu_\infty \le (|W|-1)\delta$ and claim 1 implies $\int |C_{j'}(x)(\chi_{j(w)}(wx) - C_{j(w)}(wx))| ~d\mu_\infty < \delta$. So this proves the claim.
\end{proof}
Note that the proof of claim 4 used only that $\int |\chi_a - C_a| ~d\mu_\infty \le \delta$ for all $a\in A$. So claim 3 and the proof of claim 4 imply 
$$\int |\chi''_j - C_j| ~d\mu_n \le  3|A|^2|W|\delta^{1/4} ~\forall j\in A^W.$$
Thus,
\begin{eqnarray*}
|| \phi^W_*\mu_\infty - \psi^W_*\mu_n||_1 & = & \sum_{j\in A^W} \Big| \int \chi_j ~d\mu_\infty -  \int \chi''_j ~d\mu_n\Big|\\
&\le& \sum_{j\in A^W} \int |\chi_j - C_j| ~d\mu_\infty + \Big| \int C_j ~d\mu_\infty - \int C_j ~d\mu_n\Big| + \int |\chi''_j - C_j| ~d\mu_n \\
&\le& |A^W|\big(|W|\delta + 3|A|^2|W|\delta^{1/4}\big) +\delta < \epsilon.
 \end{eqnarray*}
Since $\epsilon>0$ is arbitrary and $\Gamma \cc (X,\mu_n)$ is weakly contained in $\Gamma \cc (Y,\nu)$ this implies that $\Gamma \cc (Y,\nu)$ weakly contains $\Gamma \cc (X,\mu_\infty)$ as required.

\end{proof}

\begin{defn}
Let $(K,\kappa)$ be a standard probability space. Let $K^\Gamma$ be the set of all functions $x:\Gamma \to K$ with the product Borel structure. Let $\kappa^\Gamma$ be the product measure on $K^\Gamma$. Then $\Gamma$ acts on $K^\Gamma$ by 
$$(\gamma \cdot x)(g):=x(\gamma^{-1}g)~\forall \gamma,g\in \Gamma, x\in K^\Gamma.$$
To avoid trivialities, we assume that $\kappa$ is not concentrated on a single point.
The action $\Gamma \cc (K^\Gamma,\kappa^\Gamma)$ is called the {\em Bernoulli shift action} over $\Gamma$ with {\em base space} $(K,\kappa)$. \end{defn}

\begin{thm}
The action $\Gamma \cc (X_f,\mu_f)$ is weakly equivalent to a Bernoulli shift action.
\end{thm}
\begin{proof}
In [AWxx] it is proven that every essentially free action of any countable group $\Gamma$ weakly contains a Bernoulli shift action over $\Gamma$. In particular, all Bernoulli shift actions over $\Gamma$ are weakly equivalent. An easy exercise shows that weak containment is monotone with respect factor maps (this is also proven in [Ke09]). That is, if $\Gamma \cc (Y,\nu)$ factors onto $\Gamma \cc (Z,\zeta)$ then $\Gamma \cc (Y,\nu)$ weakly contains $\Gamma \cc (Z,\zeta)$. Thus nontrivial factors of Bernoulli shifts are weakly equivalent to Bernoulli shifts. Since $\xi:I_n^\Gamma \to X_f$ is a nontrivial factor map and $\Gamma \cc (I_n^\Gamma,\nu_n)$ is a Bernoulli shift action, it follows that $\Gamma \cc (X_f,\xi_*\nu_n)$ is weakly equivalent to a Bernoulli shift action. So the theorem follows from proposition \ref{prop:weak*} and lemma \ref{lem:Bernoulli}.
\end{proof}

\begin{defn}
Let $\Gamma \cc (Y,\nu)$. A sequence $\{A_i\}$ of Borel subsets of $Y$ is {\em asymptotically invariant} if for any $g\in \Gamma$,  $\nu(A_i \Delta gA_i) \to 0$ as $i \to\infty$. It is nontrivial if $\limsup_i \nu(A_i)(1-\nu(A_i)) > 0$. The action $\Gamma \cc (Y,\nu)$ is said to be {\em strongly ergodic} if there does not exist any nontrivial asymptotically invariant sequences. It is easy to see that strong ergodicity implies ergodicity. By [Sc81], if $\Gamma$ is amenable then no action of $\Gamma$ is strongly ergodic. 
\end{defn}

\begin{cor}\label{cor:stronglyergodic}
If $\Gamma$ is non-amenable then $\Gamma \cc (X_f,\mu_f)$ is strongly ergodic.
\end{cor}
\begin{proof}
By [LR81], Bernoulli shifts over $\Gamma$ are strongly ergodic (see also [KT08]). Strong ergodicity is preserved by weak equivalence, so it follows from the previous theorem. 
\end{proof}

\section{Entropy and periodic points: lower bound }\label{sec:lowerbound}
         The purpose of this section is to prove theorem \ref{thm:lowerbound} below and use this to finish the proofs of theorems \ref{thm:periodic} and \ref{thm:main}.
   
\begin{thm}\label{thm:lowerbound}
   If $\{\Gamma_i\}_{i=1}^\infty$ is a  sequence of finite-index normal subgroups of $\Gamma$ with $\lim_{i\to\infty}\Gamma_i=\{e\}$ and $\Gamma$ is non-amenable then $\gr\big(\{\Gamma_i\}_{i=1}^\infty,X_f\big) \le h(\{\Gamma_i\}_{i=1}^\infty,X_f,\mu_f)$.
\end{thm}                              
                    
First, we need the following lemma.

\begin{lem}
Let $\rho_n$ be the uniform probability measure on $\per(\Gamma_n,X_f)$. Then $\rho_n$ converges to $\mu_f$ in the weak* topology.
\end{lem}

\begin{proof}
Let $x\in X_f$. We claim that there exist elements $y_i \in \per(\Gamma_i,X_f)$ such that $\lim_i y_i =x$. Since $\bigcap_{n \ge 1} \bigcup_{i=n}^\infty \Gamma_i = \{e\}$, there exists an increasing sequence $\{W_i\}_{i=1}^\infty$ of finite subsets $W_i\subset \Gamma$ such that $\cup_{i=1}^\infty W_i = \Gamma$ and $W_i \cap \Gamma_i = \{e\}$.  Define $z_i \in I_M^\Gamma$ by $z_i(\gamma w)=P(x)(w)$ for all $w\in W_i$ and $\gamma \in \Gamma_i$ while $z_i(g)=0$ if $g \notin \Gamma_i W_i$. Since $W_i \nearrow \Gamma$, $\lim_{i\to\infty} z_{i} = P(x)$. Since $\xi$ is continuous, $\lim_{i\to\infty} \xi(z_{i}) = \xi(P(x))=x$. However, $y_i:=\xi(z_{i}) \in \per(\Gamma_i,X_f)$ by construction. This proves the claim.

Define $A_i:X_f \to X_f$ by $A_i(z)=z+y_i$. Since $y_i \in \per(\Gamma_i,X_f)$, $(A_i)_*\rho_i = \rho_i$. Let $\omega$ be a weak* limit point of $\{\rho_i\}_{i=1}^\infty$. Let $A_x:X_f \to X_f$ be the map $A_x(z)=x+z$. Since $y_i$ converges to $x$, $(A_x)_*\omega = \omega$. So, $\omega$ is invariant under addition by elements of $X_f$. By the uniqueness of Haar measure, this implies $\omega=\mu_f$.
\end{proof}

\begin{proof}[Proof of theorem \ref{thm:lowerbound}]
By theorem \ref{thm:alternative} and corollary \ref{cor:stronglyergodic}, it suffices to prove that $\gr\big(\{\Gamma_i\}_{i=1}^\infty,X_f\big) \le \bh(\{\Gamma_i\}_{i=1}^\infty,\phi)$. As in the previous lemma, let $\rho_n$ be the uniform probability measure on $\per(\Gamma_n)$. Because each $\Gamma_i$ is normal, $L_*(\rho_i)$ is a $\Gamma$-invariant measure supported on $[-1+\kappa,\kappa]^\Gamma$. Let $\omega$ be any weak* limit point of $\{L_*(\rho_i)\}_{i=1}^\infty$. 

\noindent {\bf Claim}. $\omega = L_*\mu_f$. 
\begin{proof}
Since the projection map $\pi^\Gamma:\R^\Gamma \to \T^\Gamma$ is continuous, $\pi^\Gamma_*\omega$ is a weak* limit point of $\{(\pi^\Gamma L)_*\rho_i\}_{i=1}^\infty$. The latter equals $\{\rho_i\}_{i=1}^\infty$ by definition. So the previous lemma implies $\pi^\Gamma_*\omega = \mu_f$. As noted in the proof of theorem \ref{thm:upperbound},
$$\mu_f\big(\{ x \in \T^\Gamma:~x(g)\ne \kappa ~\mod 1~ \forall g\in \Gamma\}\big) =1.$$
Hence, $\omega\big((-1+\kappa,\kappa)^\Gamma\big) =1.$ Since $L \circ \pi^\Gamma$ is the identity map on the set $(-1+\kappa,\kappa)^\Gamma$, it follows that $\omega=L_*\mu_f$ as claimed.
 \end{proof}
Since $P$ is the composition of $L$ with a continuous map, it follows that $\{P_*\rho_i\}_{i=1}^\infty$ converges to $P_*\mu_f$ in the weak* topology. Hence
$$ \bh(\{\Gamma_i\}_{i=1}^\infty,\phi) \ge h\big( \{ P_*\rho_i\} \big)= \limsup_n \frac{H(\rho_n)}{[\Gamma:\Gamma_n]}.$$
Since $\rho_n$ is uniformly distributed on $\per(\Gamma_n,X_f)$, $H(\rho_n) = \log |\per(\Gamma_n,X_f)|$. This implies the theorem.
\end{proof}

Theorem \ref{thm:periodic} follows immediately from theorems \ref{thm:upperbound} and \ref{thm:lowerbound}. Theorem \ref{thm:main} follows from theorem \ref{thm:periodic}, [DS07, corollary 5.3] and [DS07, theorem 6.1]. The latter two results imply that $\gr(\{\Gamma_i\}_{i=1}^\infty,X_f) = \log \det (f)$.

\end{document}